\documentclass[letterpaper, 10 pt, conference]{preamble/ieeeconf}  

\IEEEoverridecommandlockouts                              
\title{\LARGE \bf
ESDIRK-based nonlinear model predictive control \\ for stochastic differential-algebraic equations}

\author{Anders Hilmar Damm Christensen, Nicola Cantisani, John Bagterp Jørgensen
\thanks{A. H. D. Christensen, N. Cantisani, and J. B. Jørgensen are with the Department of Applied
Mathematics and Computer Science, Technical University of Denmark, DK-2800 Kgs. Lyngby, Denmark.
Corresponding author: J. B. Jørgensen (E-mail: jbjo@dtu.dk).
        {\tt\small jbjo@dtu.dk}}%
}

\usepackage{graphicx}     
\usepackage{amsmath}
\usepackage{amsfonts}

\begin{document}

\maketitle
\thispagestyle{empty}
\pagestyle{empty}


\begin{abstract}

In this paper, we present a nonlinear model predictive control (NMPC) algorithm for systems modeled by semi-explicit stochastic differential-algebraic equations (DAEs) of index 1. The NMPC combines a continuous-discrete extended Kalman filter (CD-EKF) with an optimal control problem (OCP) for setpoint tracking. We discretize the OCP using direct multiple shooting. We apply an explicit singly diagonal implicit Runge-Kutta (ESDIRK) integration scheme to solve systems of DAEs, both for the one-step prediction in the CD-EKF and in each shooting interval of the discretized OCP. The ESDIRK methods use an iterated internal numerical differentiation approach for precise sensitivity computations. These sensitivities are used to provide accurate gradient information in the OCP and to efficiently integrate the covariance information in the CD-EKF.
Subsequently, we present a simulation case study where we apply the NMPC to a simple alkaline electrolyzer stack model. We use the NMPC to track a time-varying setpoint for the stack temperature subject to input bound constraints.

\end{abstract}
\section{Introduction}
\label{sec:Introduction}

\noindent
Nonlinear model predictive control (NMPC) algorithms are widely used in the control of chemical systems to optimize economics or minimize deviations from setpoints \cite{Ritschel:Joergensen:2019}. An NMPC algorithm is a feedback control strategy that applies the solutions from optimal control problems (OCPs) based on a mathematical model of the system in a moving horizon fashion. Since chemical systems are best modeled using differential-algebraic equations (DAEs) \cite{Biegler:etal:2012}, it follows that the models employed in NMPC strategies should also be formulated as DAEs. 

Implementing the OCPs using direct multiple shooting discretization requires numerical integration schemes to solve initial-value problems (IVPs). Software solutions for solving IVPs for DAEs include the codes DASSL, which utilizes a backward differentiation formula (BDF), and RADAU5, employing a 3-stages fully implicit Runge-Kutta (IRK) integration scheme \cite{Ascher:Petzold:1998}. Due to the frequent discontinuities inherent in multiple shooting-based NMPC, BDF methods may experience frequent restarts, leading to a reduction in integration order \cite{Kristensen:etal:2004}. In contrast, the fully implicit IRK methods do not suffer from the reduction of integration order. However, these methods are in general considered computationally expensive methods as they require the solution of large systems of nonlinear equations in each integration step. Therefore, the diagonal IRK methods, such as explicit singly diagonal implicit Runge-Kutta (ESDIRK) integration schemes, offer a computationally efficient alternative. ESDIRK methods are characterized by their identical coefficients in the diagonal of their Butcher tableaus with the first step being explicit. Such methods enable the reuse of matrix factorizations across stages, thereby enhancing computational efficiency in solving IVPs for DAEs. 

Solving the OCPs also requires integrator sensitivity information. Equipping the implementation of the ESDIRK integration schemes with iterated sensitivity computations based on the principle of internal numerical differentiation (IND), ensures precise sensitivity information. The principle
of IND involves computing the sensitivities by directly differentiating the discretization scheme generated adaptively
from the integrator \cite{Bock:1981, Albersmeyer:2010}. \cite{Christensen:Joergensen:2024} presents iterated IND sensitivity calculations for ESDIRK-based OCPs applied to ODEs. It is demonstrated that such sensitivities can enhance convergence properties in an NLP solver compared to approximate sensitivities.\\
\indent
Extended Kalman filters are widely used state estimators in NMPC applications \cite{Joergensen:ecce:2007}.
\cite{Joergensen:inbook:2007} presents a continuous-discrete extended Kalman filter (CD-EKF) for systems described by stochastic DAEs (SDAEs). This CD-EKF implementation is computationally efficient due to the use of ESDIRK methods for simultaneous integration of the mean-covariance pair. The implementation obtains the covariance information by using sensitivity computations of the ESDIRK method. Here, the sensitivities are computed using staggered direct approach based on the assumption of a constant Jacobian in each accepted integration step.\\

\noindent  
This paper presents an ESDIRK-based NMPC algorithm for setpoint tracking in continuous-discrete semi-explicit index-1 SDAEs. The NMPC algorithm combines a CD-EKF with an OCP discretized using direct multiple shooting. The novelty in this paper is the direct demonstration of how the ESDIRK methods with iterated IND sensitivity computations can be applied to address both the state estimation problem and the OCP, resulting in a computationally efficient NMPC solution for SDAEs. We test the NMPC algorithm on a numerical case study of an alkaline electrolyzer stack affected by a stochastic process. We use the NMPC to track a time-varying setpoint for the temperature in the stack. 

The rest of the paper is structured as follows. Section \ref{sec:Model} introduces the SDAEs. In Section \ref{sec:simulation_and_sensitivity_gen}, we describe the ESDIRK integration scheme for solving the IVPs of DAEs in the CD-EKF and the OCP. We derive iterated IND sensitivities for these ESDIRK methods. We also introduce an implicit-explicit scheme for solving the IVPs for SDAE systems. In Section \ref{sec:NMPC}, we introduce the NMPC and demonstrate how to implement the CD-EKF and OCP for index-1 semi-explicit DAEs. Section \ref{sec:case_study} presents a case study using the model of an alkaline electrolyzer stack. Section \ref{sec:conclusions} concludes the paper.
\section{Stochastic differential-algebraic equations}
\label{sec:Model}

\noindent
We consider the continuous-discrete semi-explicit SDAEs
of the form 
\begin{subequations}
\label{eq:SDAE_model_w_meas_and_output_eqs}
\begin{alignat}{2}
    dx(t) &= f(t, x(t), y(t), u(t), d(t))dt + \sigma d\omega(t), \label{subeq:SDE}\\
    0 &= g(t, x(t), y(t), u(t), d(t)),\\
    y_k^m &= m(t_k, x_k, y_k, u_k, d_k) + v_k,\label{subeq:meas_eq}\\
    z(t) &=  h(t, x(t), y(t), u(t), d(t)) \label{subeq:out_eq}.
\end{alignat}
\end{subequations}
$t$ is time, $x(t)$ is the vector of differential states, $y(t)$ is the vector of algebraic states, $u(t)$ is the vector of inputs, and $d(t)$ is the vector of disturbances.
$\omega(t)$ is a standard Wiener process, i.e., $d\omega(t) \sim N_{iid}(0, Idt)$, representing the plant-model mismatch and $v_k \sim N_{iid}(0, R_k)$ is the measurement noise with $R_k$ being the covariance. $y_k^m = y^m(t_k)$ represent the measurements at the time $t_k$ and $z(t)$ is the output of the model. We assume the system in \eqref{eq:SDAE_model_w_meas_and_output_eqs} is of index-1, i.e., $(\partial g/ \partial y)^{-1}$ exists. 
\section{Integration and sensitivity generation}
\label{sec:simulation_and_sensitivity_gen}

\noindent
In this section, we introduce the ESDIRK integration schemes to solve IVPs for DAEs of index-1. For these ESDIRK methods, we demonstrate how to generate the forward sensitivities using the principle of IND. 
\begin{table}[tb]
\addtolength{\tabcolsep}{-3pt}    
\caption{Butcher tableaus for some ESDIRK methods.}
\label{tab:butcher_tab_ESDIRK}
    \begin{minipage}[t]{.21\linewidth}
        \vspace{0pt} 
        \centering
        
        \begin{tabular}{ccc}
        \multicolumn{3}{c}{ESDIRK12}       \\ \\
        \multicolumn{1}{c|}{$0$} & $0$ &   \\
        \multicolumn{1}{c|}{$1$} & $b_1$ & $\gamma$ \\ \hline
        \multicolumn{1}{c|}{$x_{k+1}$} & $b_1$ & $\gamma$ \\ \hline
        \multicolumn{1}{c|}{$\hat{x}_{k+1}$} & $\hat{b}_1$ & $\hat{b}_2$ 
        \end{tabular}
    \end{minipage}
    \begin{minipage}[t]{.33\linewidth}
        \vspace{0pt} 
        \centering
        
        \begin{tabular}{cccc}
        \multicolumn{4}{c}{ESDIRK23}       \\ \\
        \multicolumn{1}{c|}{$0$} & $0$ &    \\
        \multicolumn{1}{c|}{$c_2$} & $a_{21}$ & $\gamma$ &  \\ 
        \multicolumn{1}{c|}{$1$} & $b_1$ & $b_1$ & $\gamma$ \\ \hline
        \multicolumn{1}{c|}{$x_{k+1}$} & $b_1$ & $b_1$ & $\gamma$ \\ \hline
        \multicolumn{1}{c|}{$\hat{x}_{k+1}$} & $\hat{b}_1$ & $\hat{b}_2$ & $\hat{b}_3$
        \end{tabular}
    \end{minipage}
    \begin{minipage}[t]{.35\linewidth}
        
        \vspace{0pt} 
        \centering
        
        \begin{tabular}{ccccc}
        \multicolumn{5}{c}{ESDIRK34}       \\ \\
        \multicolumn{1}{c|}{$0$} & $0$ &   \\
        \multicolumn{1}{c|}{$c_2$} & $a_{21}$ &  $\gamma$ &  \\ 
        \multicolumn{1}{c|}{$c_3$} & $a_{31}$ & $a_{32}$ & $\gamma$ &  \\ 
        \multicolumn{1}{c|}{$1$} & $b_1$ & $b_2$ & $b_3$ & $\gamma$\\ \hline
        \multicolumn{1}{c|}{$x_{k+1}$} & $b_1$ & $b_2$ & $b_3$ & $\gamma$ \\ \hline
        \multicolumn{1}{c|}{$\hat{x}_{k+1}$} & $\hat{b}_1$ & $\hat{b}_2$ & $\hat{b}_3$ & $\hat{b}_4$ 
        \end{tabular}
    \end{minipage}
    \addtolength{\tabcolsep}{3pt}
\end{table}

\subsection{The ESDIRK integration schemes for DAEs}
\noindent
We solve the IVPs
\begin{subequations}
\label{eq:DAEs}
\begin{alignat}{2}
    \frac{dx(t)}{dt} &= f(t, x(t), y(t), u, d), \quad x(t_0) = x_0,\label{subeq:ODE}\\
    0 &= g(t, x(t), y(t), u, d), \quad y(t_0) = y_0,\label{subeq:Algebraic}
\end{alignat}
\end{subequations}
for $t\in [t_0,\, t_f]$ using the ESDIRK integration schemes. The $n_s$-stages ESDIRK methods for solving \eqref{subeq:ODE} are \cite{Kristensen:etal:2005}
\begin{subequations}
    \label{eq:rk_scheme}
\begin{alignat}{2}
    T_i &= t_k+c_ih_k,\label{subeq:nodes}\\
    X_i &= x_k+h_k\sum_{j=1}^{i}a_{ij}f(T_j, X_j, Y_j, u, d), \label{subeq:stages}\\
    x_{k+1} &= x_k+h_k\sum_{i=1}^{n_s}b_if(T_j, X_j, Y_j, u, d),\label{subeq:nextIter}\\
    \hat{x}_{k+1} &= x_k+ h_k\sum_{i=1}^{n_s}\hat{b}_if(T_j, X_j, Y_j, u, d). \label{subeq:bhat}
\end{alignat}
\end{subequations}
$T_i$ are the internal nodes, $X_i$ are the differential stages, and $Y_i$ algebraic stages at iteration $k$ for $i = 1, \dots, n_s$. $x_k$ and $x_{k+1}$ are the steps computed at $t_k$ and $t_{k+1} = t_k + h_k$. $h_k$ is the integration step size at iteration $k$ which may be controlled based on a local error estimate, $e_{k+1} = x_{k+1}-\hat{x}_{k+1}$, using the embedded method. To satisfy the consistency equations in  \eqref{subeq:Algebraic} we compute $Y_j$ such that
\begin{equation}
\label{eq:algebraic_rk_scheme}
    0 = g(T_j, X_j, Y_j, u, d), \quad \text{for} \quad j=1,\dots, n_s.
\end{equation}
Table \ref{tab:butcher_tab_ESDIRK} presents the Butcher tableaus for the ESDIRK12, ESDIRK23, and ESDIRK34 methods. The numerical values are provided by \cite{Joergensen:arxiv:2018}. For the rest of this section, we leave out the arguments $T_j$ and $d$.

We combine \eqref{subeq:stages} and \eqref{eq:algebraic_rk_scheme} into the system of nonlinear equations
\begin{equation}
\label{eq:ESDIRK_residual}
R_i(S_i) =\begin{bmatrix}
 X_i-h_k\gamma f(X_i, Y_i, u)-\psi_i \\ 
- g(X_i, Y_i, u)
\end{bmatrix} = 0, 
\end{equation}
for $i = 2, \dots, n_s$. $S_i = \begin{bmatrix}
 X_i; & Y_i
\end{bmatrix}$ and
\begin{equation}
\psi_i = x_k+h_k\sum_{j=1}^{i-1}a_{ij}f(X_j, Y_j, u),
\end{equation}
and solve \eqref{eq:ESDIRK_residual} in each integration step. We do this using the inexact Newton scheme
\begin{subequations}
    \label{eq:newton_method}
\begin{alignat}{2}
    M_k\Delta S_i^{[l]} &= -R_i( S_i^{[l]}),\\
    S_i^{[l+1]} &= S_i^{[l]} +\Delta S_i^{[l]},
\end{alignat}
\end{subequations}
for $l = 0, \dots, v_{i,k}-1$. $v_{i,k}$ is the number of Newton-type iterations required for stage $i$ at iteration $k$ to satisfy some convergence criteria. We choose the convergence criteria
\begin{equation}
\label{eq:stop_crit_newton}
    \left \| R_i(S_i^{[l]}) \right \| = \max_{j\in 1, \dots , (n_S)}\frac{|R_i(S_i^{[l]})_j|}{\max(\text{abs}, \,\text{rel}(S_i^{[l]})_j)} < \tau,
\end{equation}
where $n_S$ is the dimension of $S_i$ and we choose $\tau = 0.1$.
$M_k \approx \frac{\partial R_i(S_i^{[l]})}{\partial S_i^{[l]}}$ is the iteration matrix defined as
\begin{equation}
\label{eq:iteration_matrix}
M_k = \begin{bmatrix}
 I-h_k\gamma \frac{\partial f(x_k, y_k, u) }{\partial x_k}& -h_k\gamma \frac{\partial f(x_k, y_k, u) }{\partial y_k}\\ 
 -\frac{\partial g(x_k, y_k, u) }{\partial x_k} & -\frac{\partial g(x_k, y_k, u) }{\partial y_k}
\end{bmatrix},
\end{equation}
where $I$ is the identity matrix. For fixed step size ESDIRK methods, we update $M_k$ at every integration step and we use its LU factorizations to solve for all the stages in \eqref{eq:newton_method}.
Since the ESDIRK integration schemes are stiffly accurate Runge-Kutta methods, we directly obtain the next differential and algebraic states as 
\begin{equation}
\label{eq:stiffly_accurate_next_iter}
    x_{k+1}=X_{n_s}^{v_{n_s,k}-1}, \quad y_{k+1}=Y_{n_s}^{v_{n_s,k}-1}.
\end{equation}

\noindent
We apply stage value predictors (SVPs) of the form
\begin{equation}
\label{eq:SVPs}
    S_i^{[0]} = \alpha_i(r_k)s_{k-1} + \sum_{j=2}^{n_s}\beta(r_k)_{ij}\hat{S_j},
\end{equation}
to provide initial guesses for the Newton-type iterations in \eqref{eq:newton_method}. $\alpha(r_k) \in \mathbb{R}^{n_s-1}$ and $\beta(r_k) \in \mathbb{R}^{n_s-1\times n_s-1}$ are the predictor coefficients computed using $r_k = h_{k}/h_{k-1}$. 
$\hat{S_j} = \hat{S_j}^{[v_{j,k-1}-1]} = \begin{bmatrix}
\hat{X_j}^{[v_{j,k-1}-1]}; & \hat{Y_j}^{[v_{j,k-1}-1]}
\end{bmatrix}$ for $j=2, \dots n_s$ represent the previously converged differential and algebraic stages and $s_{k-1} = \begin{bmatrix}
x_{k-1}; & y_{k-1}
\end{bmatrix}$. We apply the order conditions in \cite{Higueras:Roldan:2005} to compute the predictor coefficients  $\alpha(r_k)$ and $\beta(r_k)$, and we set $r_k = 1$, as we only consider fixed step-size ESDIRK methods.
\cite{Christensen:Joergensen:2024} shows the predictor coefficients for ESDIRK12 and ESDIRK23. For the first step, $k=0$, we apply the trivial predictor $S^{[0]}_i = s_k = \begin{bmatrix}
x_{k}; & y_{k}
\end{bmatrix}$.

\subsection{Iterated IND sensitivity generation}

\noindent
We differentiate the adaptively generated discretization scheme of the ESDIRK methods. For fixed step-size ESDIRK methods, the adaptive components are the LU factorizations of the iteration matrices, $M_k$, the number of Newton-type iterations for all stages, $v_{i,k}$, and the sequence of Newton-type iterates, $S_i^{[l]}$ for $l = 0, \dots, v_{i,k}-1$ and for $i = 2, \dots, n_s$.
We define the sensitivities of the combined differential and algebraic states with respect to initial states and inputs as
\begin{equation}
\label{eq:sensitivities_k}
    \frac{\partial s_{k}}{\partial s_0} = \begin{bmatrix}
\frac{\partial x_k}{\partial x_0} & \frac{\partial x_k}{\partial y_0}\\ 
 \frac{\partial y_k}{\partial x_0}&  \frac{\partial y_k}{\partial y_0}
\end{bmatrix}, \quad \frac{\partial s_{k}}{\partial u} = \begin{bmatrix}
\frac{\partial x_k}{\partial u}\\ 
\frac{\partial y_k}{\partial u}
\end{bmatrix},
\end{equation}
with initial conditions $\frac{\partial s_{0}}{\partial s_0} = I$ and $\frac{\partial s_{0}}{\partial u} = 0$. We express the implementation of the ESDIRK integration schemes in terms of the elementary operations
\begin{subequations}
    \label{eq:esdirk_elementary_operations}
\begin{alignat}{2}
    S_i^{[0]} &= \phi^{\mathrm{SVP}}_i(s_{k-1}, \hat{S}), \label{subeq:initGuess}\\
    S_i^{[l+1]} &= \phi^{\mathrm{it}}_i( S_i^{[l]}, \psi_i, u) \nonumber\\
    &= S_i^{[l]} - M_k^{-1}R_i(S_i^{[l]}, \psi_i, u),\label{subeq:NIter}\\
    s_{k+1} &= \phi^{n_s}(S_{n_s}^{[v_{n_s,k}-1]}) = S_{n_s}^{[v_{n_s,k}-1]}\label{subeq:stage2nextstate},
\end{alignat}
\end{subequations}
where \eqref{subeq:initGuess} represents the SVPs for stage $i$ with \\
$\hat{S}~=~\begin{bmatrix}
\hat{S}_2^{[v_{2,k-1}-1]}; & \dots; & s_k
\end{bmatrix}$. Eq. \eqref{subeq:NIter} is the Newton-type scheme in \eqref{eq:newton_method} and \eqref{subeq:stage2nextstate} represent the \eqref{eq:stiffly_accurate_next_iter}.
\subsubsection{State sensitivity}
We construct the state sensitivities by computing the partial derivatives of \eqref{eq:esdirk_elementary_operations} with respect to the initial differential and algebraic states, $s_0$.
The sensitivities related to \eqref{subeq:initGuess} are
\begin{equation}
    \frac{\partial S_i^{[0]}}{\partial s_0} = \alpha_i(r_k)\frac{\partial s_{k-1}}{\partial s_0} + \sum_{j=2}^{n_s} \beta_{ij}(r_k)\frac{\partial \hat{S}_j}{\partial s_0},
\end{equation}
for $k \geq 1$ and $\frac{\partial S_i^{[0]}}{\partial s_0} = \frac{\partial s_k}{\partial s_0}$ for $k=0$. For \eqref{subeq:NIter}, we determine the state sensitivities using the Newton-type scheme
\begin{equation}
    \frac{\partial S_i^{[l+1]}}{\partial s_0} = \frac{\partial S_i^{[l]}}{\partial s_0}-M_k^{-1}\frac{\partial R_i(S_i^{[l]}, \psi_i, u)}{\partial s_0},
\end{equation}
where $\frac{\partial R_i(S_i^{[l]}, \psi_i, u)}{\partial s_0}$ is
\begin{equation}
\label{eq:dRdz0}
    \frac{\partial R_i(S_i^{[l]}, \psi_i, u)}{\partial s_0} = \frac{\partial R_i}{\partial S_i^{[l]}}\frac{\partial S_i^{[l]}}{\partial s_0} + \frac{\partial R_i}{\partial \psi_i}\frac{\partial \psi_i}{\partial s_0}.
\end{equation}

The terms in \eqref{eq:dRdz0} are 
\begin{subequations}
\label{eq:partial_derivatives_in_ESDIRK_methods}
\begin{alignat}{2}
    \frac{\partial R_i}{\partial S_i^{[l]}} &= \begin{bmatrix}
 I-h_k\gamma \frac{\partial f(X_i^{[l]},Y_i^{[l]}, u) }{\partial X_i^{[l]}}& -h_k\gamma \frac{\partial f(X_i^{[l]},Y_i^{[l]}, u) }{\partial Y_i^{[l]}}\\ 
 -\frac{\partial g(X_i^{[l]},Y_i^{[l]}, u) }{\partial X_i^{[l]}} & -\frac{\partial g(X_i^{[l]},Y_i^{[l]}, u) }{\partial Y_i^{[l]}}
\end{bmatrix}, \label{subeq: dRdXi}\\
    \frac{\partial R_i}{\partial \psi_i} &= \begin{bmatrix} -I \\ 0\end{bmatrix}\label{subeq: dRdpsi}, \quad
    \frac{\partial \psi_i}{\partial s_0} = \begin{bmatrix} \frac{\partial \psi_i}{\partial x_0} & \frac{\partial \psi_i}{\partial y_0}\end{bmatrix}.
\end{alignat}
\end{subequations}
The partial derivatives in $\frac{\partial \psi_i}{\partial s_0}$ uses the Jacobians of the previously converged stages for iteration $k$ as
\begin{subequations}
\begin{alignat}{2}
   \frac{\partial \psi_i}{\partial x_0} &= \frac{\partial x_k}{\partial x_0}+h_k\sum_{j=1}^{i-1}a_{ij}\bigg(\frac{\partial f(X_j, Y_j, u)}{\partial S_j}\frac{\partial S_j}{\partial x_0}\bigg)\\
   \frac{\partial \psi_i}{\partial y_0} &= \frac{\partial x_k}{\partial y_0}+h_k\sum_{j=1}^{i-1}a_{ij}\bigg(\frac{\partial f(X_j, Y_j, u)}{\partial S_j}\frac{\partial S_j}{\partial y_0}\bigg),
\end{alignat}
\end{subequations}
where
\begin{equation}
    \frac{\partial f(X_j, Y_j, u)}{\partial S_j} = \begin{bmatrix}
\frac{\partial f(X_j, Y_j, u)}{\partial X_j} &  \frac{\partial f(X_j, Y_j, u)}{\partial Y_j}
\end{bmatrix}.
\end{equation}

\noindent
Finally, due to the ESDIRK methods being stiffly accurate, we directly obtain the sensitivities for \eqref{subeq:stage2nextstate}
as
\begin{equation}
    \frac{\partial x_{k+1}}{\partial x_0} = \frac{\partial X_{n_s}^{v_{n_s},k-1}}{\partial x_0}, \quad  \frac{\partial x_{k+1}}{\partial y_0} = \frac{\partial X_{n_s}^{v_{n_s},k-1}}{\partial y_0}.
\end{equation}

\subsubsection{Input sensitivity}
Similarly to the state sensitivities, we generate the input sensitivities by computing the partial derivatives of \eqref{eq:esdirk_elementary_operations} with respect to the vector of inputs. For the initial guess of the Newton-type scheme we derive the sensitivities at $k \geq 1$ as
\begin{equation}
    \frac{\partial S_i^{[0]}}{\partial u} = \alpha_i(r_k)\frac{\partial s_{k-1}}{\partial u} + \sum_{j=2}^{n_s} \beta_{ij}(r_k)\frac{\partial \hat{S}_j}{\partial u},
\end{equation}
with $\frac{\partial S_i^{[0]}}{\partial u} = \frac{\partial s_k}{\partial u}$ for $k = 0$. For the iterations in \eqref{subeq:NIter} we compute the input sensitivity
\begin{equation}
    \frac{\partial S_i^{[l+1]}}{\partial u} = \frac{\partial S_i^{[l]}}{\partial u}-M_k^{-1}\frac{\partial R_i(S_i^{[l]}, \psi_i, u)}{\partial u},
\end{equation}

where the partial derivative of the residual function is
\begin{equation}
    \frac{\partial R_i(S_i^{[l]}, \psi_i, u)}{\partial u} = \frac{\partial R_i}{\partial S_i^{[l]}}\frac{\partial S_i^{[l]}}{\partial u} + \frac{\partial R_i}{\partial \psi_i}\frac{\partial \psi_i}{\partial u} + \frac{\partial R_i}{\partial u}.
\end{equation}

\noindent
The expressions $ \frac{\partial R_i}{\partial S_i^{[l]}}$ and $\frac{\partial R_i}{\partial \psi_i}$ are shown in \eqref{eq:partial_derivatives_in_ESDIRK_methods} while we calculate the other terms as

\begin{subequations}
\begin{alignat}{2}
   \frac{\partial \psi_i}{\partial u} &= \frac{\partial x_k}{\partial u} \nonumber\\
   &+h_k\sum_{j=1}^{i-1}a_{ij}\bigg(\frac{\partial f(X_j, Y_j, u)}{\partial S_j}\frac{\partial S_j}{\partial u}+ \frac{\partial f(X_j, Y_j, u)}{\partial u}\bigg), \\
   \frac{\partial R_i}{\partial u}  &=-\begin{bmatrix}
h_k\gamma\frac{\partial f(X_i^{[l]},Y_i^{[l]}, u) }{\partial u}; & \frac{\partial g(X_i^{[l]},Y_i^{[l]}, u) }{\partial u}
\end{bmatrix}.
\end{alignat}
\end{subequations}

We directly obtain the input sensitivities for the next state as $\frac{\partial x_{k+1}}{\partial u} = \frac{\partial X_{n_s}^{v_{n_s},k-1}}{\partial u}$.

\subsection{Implicit-explicit SDAE solver}
We solve \eqref{eq:SDAE_model_w_meas_and_output_eqs} using an implicit-explicit integration scheme \cite{Ritschel:Joergensen:2019}. 
The implicit-explicit integration scheme applies $M$-steps in each interval $t\in [t_k, t_{k+1}]$ for $k = 0, \dots, N-1$ and is described by
\begin{subequations}
    \label{eq:implicit_explicit_DAE}
\begin{alignat}{2}
    t_{k,n+1} &= t_{k,n} + \Delta t_{k,n},\\
    x_{k,n+1} &= x_{k,n} + f(t_{k,n+1}, x_{k,n+1}, y_{k,n+1}, u_k, d_k)\Delta t_{k,n} \nonumber\\
    &+ \sigma \Delta \omega_{k,n},\\
    0 &= g(t_{k,n+1}, x_{k,n+1}, y_{k,n+1}, u_k, d_k),
\end{alignat}
\end{subequations}
for $n = 0, \dots, M-1$ and where $t_{k,0}=t_k$, $x_{k,0}=x_k$, $y_{k,0}=y_k$, and $\Delta\omega_{k,n} \sim N_{iid}(0, I\Delta t_{k,n})$. The integration step-size is $\Delta t_{k,n} = t_{k,n+1}-t_{k,n}$ and  we obtain the step at $t_{k+1}$ as
\begin{equation}
    t_{k+1} = t_{k,M}, \quad  x_{k+1} = x_{k,M},\quad  y_{k+1} = y_{k,M}.
\end{equation}

\noindent
In each iteration $n$ in \eqref{eq:implicit_explicit_DAE}, we solve the system of nonlinear equations, $R(s_{k,n+1}) = 0$. We define $R(s_{k,n+1})$ as
\begin{equation}
\label{eq:residual_SDAE}
    R(s_{k,n+1}) = \begin{bmatrix}
x_{k,n+1} - x_{k,n} - f_{k,n}\Delta t_{k,n}-\sigma\omega_{k,n}\\ 
g_{k,n+1}
\end{bmatrix},
\end{equation}
where $s_{k,n+1} = \begin{bmatrix}
x_{k,n+1}; & y_{k,n+1}
\end{bmatrix}$ and
\begin{subequations}
\label{eq:short_notation_SDAE_solver}
\begin{alignat}{2}
    f_{k,n} &= f(t_{k,n+1}, x_{k,n+1}, y_{k,n+1}, u_k, d_k),\\
    g_{k,n+1} &= g(t_{k,n+1}, x_{k,n+1}, y_{k,n+1}, u_k, d_k).
\end{alignat}
\end{subequations}

We solve \eqref{eq:residual_SDAE} using an exact Newton scheme 
\begin{subequations}
    \label{eq:newton_method_SDAE}
\begin{alignat}{2}
    \frac{\partial R(s_{k,n+1}^{[l]})}{\partial s_{k,n+1}^{[l]}}\Delta s_{k, n+1}^{[l]} &= -R( s_{k,n+1}^{[l]}),\\
    s_{k,n+1}^{[l+1]} &= s_{k,n+1}^{[l]} +\Delta s_{k, n+1}^{[l]},
\end{alignat}
\end{subequations}
and use a similar convergence criteria as in  \eqref{eq:stop_crit_newton}.

\section{Nonlinear model predictive control}
\label{sec:NMPC}

\noindent
This section presents the implementations of the CD-EKF and the OCP. We combine the CD-EKF and OCP into an NMPC such that we compute one control interval between each measurement. We assume the vector of disturbances is piece-wise constant between these measurements. We, therefore, apply the zero-order-hold
parametrizations 
\begin{equation}
    u(t) = u_k, \quad  d(t) = d_k \quad \text{for} \quad t\in[t_k, \, t_{k+1}],
\end{equation}
where the sampling time is $T_s = t_{k+1}-t_k$. $u_k$ represents the implemented control signal computed by the OCP at time $t_k$.

\subsection{Continuous-discrete extended Kalman filter}

We implement the CD-EKF using the ESDIRK integration scheme for DAEs as in \cite{Joergensen:inbook:2007}.
The CD-EKF computes the filtered differential states, $\hat{x}_{k|k}$, algebraic states, $\hat{y}_{k|k}$, and differential state covariance, $P_{k|k}$, at time $t_k$ based on the measurements, $y^m_k$, and the previous implemented control, $u_{k-1}$. We use the filtered states as initial conditions for the OCP. The implemented control signal $u_{k}$ is used to compute the one-step prediction at time $t_{k+1}$ in the differential states, $\hat{x}_{k+1|k}$, algebraic states, $\hat{y}_{k+1|k}$, and the differential state covariance, $P_{k+1|k}$. The initial states and covariance are 
\label{eq:CDEKF_init}
\begin{equation}
     \hat{x}_{0|-1} = x_0, \quad  \hat{y}_{0|-1} = y_0, \quad P_{0|-1} = P_0.
\end{equation}

\subsubsection{Filtering}
We compute the innovation, $e_k$, as
\begin{equation}
\label{eq:inno}
    e_k = y_k^m-\hat{y}^m_{k|k-1}, 
\end{equation}
using the estimated measurement at time $t_k$
\begin{equation}
    \hat{y}^m_{k|k-1} = m(t_k, \hat{x}_{k|k-1}, \hat{y}_{k|k-1}, u_{k-1}, d_k).
\end{equation}
The covariance of the innovation, $R_{e,k}$, is computed as 
\begin{subequations}
\label{eq:Ck}
\begin{alignat}{2}
    R_{e,k} &= C_kP_{k|k-1}C_k'+R_k,\\
    C_k &= \frac{\partial m_{k|k-1}}{\partial x} + \frac{\partial m_{k|k-1}}{\partial y}\frac{\partial \hat{y}_{k|k-1}}{\partial x}.
\end{alignat}
\end{subequations}
The sensitivity $\frac{\partial \hat{y}_{k|k-1}}{\partial x}$ is computed by solving 
\begin{equation}
    \frac{\partial g_{k|k-1}}{\partial  y}\frac{\partial y_{k|k-1}}{\partial x} = -\frac{\partial g_{k|k-1}}{ \partial x}.
\end{equation}
In these expressions, we used the notation
\begin{subequations}
\label{eq:short_notation_filtering}
\begin{alignat}{2}
    m_{k|k-1} &= m(t_k, \hat{x}_{k|k-1}, \hat{y}_{k|k-1}, u_{k-1}, d_k),\\
    g_{k|k-1} &= g(t_k, \hat{x}_{k|k-1}, \hat{y}_{k|k-1}, u_{k-1}, d_k).
\end{alignat}
\end{subequations}

We combine \eqref{eq:inno}-\eqref{eq:Ck} to compute the filtered state 
\begin{equation}
    \hat{x}_{k|k} = \hat{x}_{k|k-1} + K_ke_k,
\end{equation}
and the filtered covariance as 
\begin{equation}
    P_{k|k} = (I-K_kC_k)P_{k|k-1}(I-K_kC_k)'+K_kR_kK_k',
\end{equation}
with $K_k$ being the Kalman gain defined as
\begin{equation}
    K_k =P_{k|k-1}C_k'R_{e,k}^{-1}.
\end{equation}
Finally, to compute the filtered algebraic state, $\hat{y}_{k|k}$, we solve the system of nonlinear equations
\begin{equation}
    0 = g(t_k,\hat{x}_{k|k},   \hat{y}_{k|k}, u_{k-1}, d_k),
\end{equation}
using an exact Newton scheme.

\subsubsection{Prediction}
We compute one-step predictions, $\hat{x}_{k+1|k} = \hat{x}_k(t_{k+1})$ and $\hat{y}_{k+1|k} = \hat{y}_k(t_{k+1})$, by solving
\begin{subequations}
\label{eq:CDEKF_prediction_mean}
\begin{alignat}{2}
    \frac{d \hat{x}_k(t)}{dt} &= f(\hat{x}_k(t), \hat{y}_k(t), u_k, d_k),\\
    0 &= g(\hat{x}_k(t), \hat{y}_k(t), u_k, d_k),
\end{alignat}
\end{subequations}
for $t \in [t_k, \, t_{k+1}]$, with initial conditions $\hat{x}_k(t) = \hat{x}_{k|k}$ and $\hat{y}_k(t) = \hat{y}_{k|k}$ using an ESDIRK method. \cite{Joergensen:inbook:2007} suggest computing the predicted differential state covariance, $P_{k+1|k} = P_k(t_{k+1})$, as the solution to 
\begin{equation}
\label{eq:CDEKF_prediction_cov_integral}
    \begin{split}
     P_k(t) &= \Phi_{xx}(t,t_k)P_k(t_k)\Phi_{xx}'(t,t_k)\\
    &+ \int_{t_k}^t\Phi_{xx}(t,s)\sigma\sigma'\Phi_{xx}(t,s)ds,
    \end{split}
\end{equation}
at time $t=t_{k+1}$, with the initial condition $P_k(t_k) = P_{k|k}$. $\Phi_{xx}(t,s) = \frac{\partial \hat{x}_k(t)}{\partial \hat{x}_k(s)}$ is the differential state sensitivity which we compute using the iterated IND approach.

\subsection{Optimal control problem}
We consider an OCP that penalizes the tracking error, $z(t)-\bar{z}(t)$, and the input rate-of-movement, $\Delta u_k = u_{k}-u_{k-1}$, in a least-squares formulation. We solve the OCP at time $t_k$ with prediction horizon $T_N$ using the filtered differential state, $\hat{x}_{k|k}$, as initial conditions. We divide the control horizon, $[t_k, \, t_k+T_N]$, into $N = T_N/T_s$ equally spaced subintervals, $[t_{k+j}, \, t_{k+j+1}]$, with  $j \in \mathcal{N} = 1,\dots, N-1$. The OCP is

\begin{subequations}
\label{eq:OCP}
\begin{alignat}{4}
       &\min_{x, y, u} \quad &&\phi_k, &&&\\
       &s.t.                 &&x(t_k) = \hat{x}_{k|k}, &&&\\
       &                     &&\frac{dx(t)}{dt} = f(t, x, y, u, d), &&&t\in [t_k, \, t_k+T_N],\\
       &                     &&0 = g(t, x, y, u, d),  &&&t\in [t_k, \, t_k+T_N], \\
       &                     &&z(t) = h(t, x, y, u, d), &&&\label{subeq:output_equation}\\ 
       &                     &&u(t) = u_{k+j}, \quad j \in \mathcal{N}, &&&t\in [t_{k+j}, \, t_{k+j+1}],\\
       &                     &&d(t) = d_{k+j},\quad j \in \mathcal{N}, &&&t\in [t_{k+j}, \, t_{k+j+1}],  \label{subeq_dist_j}\\
       &                     &&u_{\min} \leq u_{k+j} \leq u_{\max}, &&&j \in \mathcal{N},
\end{alignat}
\end{subequations}
where $\phi_k = \phi_k([x(t); y(t); u(t)]_{t_k}^{t_k+T_N})$ is the objective, and $f(t, x, y, u, d) = f(t, x(t), y(t), u(t), d(t))$ and $g(t, x, y, u, d) = g(t, x(t), y(t), u(t), d(t))$ are the DAEs. We formulate the objective  $\phi_k = \phi_{z,k} +  \phi_{\Delta u,k} + \phi_{N, k}$ with the objectives function terms
\begin{subequations}
\label{eq:objectives}
\begin{alignat}{2}
       \phi_{z,k} &= \frac{1}{2}\int_{t_k}^{t_k+T_N}\left \|  z(t)-\Bar{z}(t)\right \|_{Q_{z}}^2dt, \label{eq:cont_obj}\\
       \phi_{\Delta u,k} &= \frac{1}{2}\sum_{j=0}^{N-1}\left \|  \Delta u_{k+j}\right \|_{\bar{Q}_{\Delta u}}^2,\\
       \phi_{N,k} &= \frac{1}{2}\left \|  z(t_k+T_N)-\Bar{z}(t_k+T_N)\right \|_{\bar{Q}_{z}}^2.
\end{alignat}
\end{subequations}
$\phi_{z,k}$ penalizes the tracking error $z(t)-\bar{z}(t)$, with $\bar{z}(t)$ being the setpoint and  $Q_{z}$ being the weight matrix. $\phi_{\Delta u, k}$ penalizes the input rate-of-movement with weight $\bar{Q}_{\Delta u} = Q_{\Delta u}/T_s$. Finally, $\phi_{N,k}$ penalizes the tracking error at the end of the horizon with the scaled weight matrix  $\bar{Q}_{z} = Q_{z}/T_s$.\\

\noindent
We apply the direct multiple shooting discretization scheme in \cite{Bock:etal:2007} to transcribe \eqref{eq:OCP}-\eqref{eq:objectives} into a NLP of the form
\begin{subequations}
\label{eq:NLP_OCP}
\begin{alignat}{4}
      &  \min_{w} \quad && \phi,  &&&\\
       &s.t. &&w^x_0 = \hat{x}_{k|k}, &&& \label{eq:eqcon1}\\
       & &&R(w^x_{j+1}, w^x_{j}, w^y_{j}, u_{j}, d_{j}) = 0,  \quad  &&&j \in \mathcal{N},\label{eq:eqcon2}\\ 
       & && u_{\min} \leq u_{k+j} \leq u_{\max}, \; &&&j \in \mathcal{N},
\end{alignat}
\end{subequations}
where $R(\cdot)$ is the ESDIRK integration scheme and the decision variables are 
\begin{equation}
    w = (w^x_0,\, w^y_0,\, u_0,\, w^x_{1},\, w^y_{1},\, u_{1}, \, \dots, u_{N-1},\, w^x_{N}).
\end{equation} 
The $j \in \mathcal{N}$ are relative to $t_k$ and we initialize $w^y_0$ using the filtered algebraic state, $w^y_0 = \hat{y}_{k|k}$. Since \eqref{eq:cont_obj} depends continuously on the states, we discretize it by formulating it as an IVP, which we solve simultaneously with the ESDIRK integration scheme. To allow temporary violation of the algebraic constraints in the ESDIRK integration scheme, we solve the \textit{relaxed} DAEs \cite{Bock:etal:2007} on each subinterval,
\begin{subequations}
    \label{eq:relaxed_DAE}
\begin{alignat}{2}
       \frac{dx_{j}(t)}{dt} &= f(t, x_{j}(t), y_{j}(t), u_{j}, d_{j}),\\
       0 &= g(t, x_{j}(t), y_{j}(t), u_{j}, d_{j}) \nonumber\\
       &-p_{j}(t)g(t_{j}, w^x_{j},  w^y_{j}, u_{j}, d_{j}),
\end{alignat}
\end{subequations}
with initial conditions $x_j(t_j) = w^x_{j}$ and $y_j(t_j) = w^y_{j}$. $p_j(t)$ is a scalar function such that $p_j(t_j) = 1$ and $p_j(t) > 0$ for $ t \in [t_{k+j}, \, t_{k+j+1}]$ and we choose  $p_j(t) = \exp({-\eta \frac{t-t_j}{t_{j+1}-t_j}})$
with $\eta = 1$. \\

We solve \eqref{eq:NLP_OCP} using an SQP solver with a BFGS for updating the Lagrangian Hessian approximation. We denote the equality constraints \eqref{eq:eqcon1}-\eqref{eq:eqcon2} as $b(w)$, and the gradient, $\frac{\partial b(w)}{\partial w}$, required by the SQP solver, is computed using the sensitivity information $\frac{\partial R(w)}{\partial w}$. As a result, the SQP solver requires both state and input sensitivities.
We implement the first control signal at time $t_k$, $u^*_0$, from the solution $w^*$ to the system. For the first iteration of the NMPC, $k=0$, we initialize $w$ using the initial consistent conditions, $ w_0 = (x_0, y_0, u_0, x_0, y_0, u_0, \dots, x_0)$, but for the subsequent executions of the NLP, we initialize $w$ by shifting the previously converged solution,
$w_0 = (w^{x,*}_1, w^{y,*}_1, u_1^*, \dots, w^{x,*}_N, w^{y,*}_{N-1}, u^*_{N-1}, w^{x,*}_N)$.


\section{Numerical case study}
\label{sec:case_study}

\noindent
We apply the NMPC algorithm to an alkaline electrolyzer stack model. The model we use is a simplified version of the one presented by \cite{Cantisani2023}, similar to \cite{Rizwan:etal:2021}. The system of SDAEs are
\begin{subequations}
    \begin{alignat}{1}
    dT &= \frac{1}{C_{\mathrm{p,el}}} \Big( f_{\mathrm{in}} c_{\mathrm{P,lye}} (T_{\mathrm{in}} - T) \nonumber \\ &+ n_c (U_{\mathrm{cell}} - U_{\mathrm{tn}}) I - A_s h_c (T - T_{\mathrm{amb}}) \Big)dt, \\
    dT_{\mathrm{in}} &= \sigma d\omega, \\
    0 &= U_{\mathrm{cell}} - ( U_{\mathrm{rev}} + U_{\mathrm{ohm}}(T,I) + U_{\mathrm{act}}(T,I) ), \\
    0 &= P_{\mathrm{in}} - n_c U_{\mathrm{cell}} I,
    \end{alignat}
\end{subequations}
where
\begin{subequations} 
    \begin{alignat}{1}
    U_{\mathrm{ohm}}(T,I) & = (r_1 + r_2 T)\frac{I}{A}, \\
    U_{\mathrm{act}}(T,I) & = s \log \left( \left(t_1 + \frac{t_2}{T} + \frac{t_3}{T^2} \right) \frac{I}{A} + 1 \right).
    \end{alignat}
\end{subequations}
The SDAEs model the energy balance and electrochemical dynamics in the electrolyzer stack. $x(t) = \begin{bmatrix}
T; & T_{in}
\end{bmatrix}$  being the vector of differential variables representing the temperature in the electrolyzer stack and the inlet temperature, respectively. In this system, the inlet temperature is modeled as a stochastic differential equation with a zero drift and where $\sigma = 0.03$. $y(t) = \begin{bmatrix}
 U_{\mathrm{cell}};& I
\end{bmatrix}$ is the vector of algebraic equations with $U_{\mathrm{cell}}$ being the cell voltage and $I$ being the current. $C_{\mathrm{p,el}}$, $c_{\mathrm{P,lye}}$, $n_c,U_{\mathrm{tn}}$, $A_s$, $h_c$, $U_{\mathrm{rev}}$, and $A$ are parameters of the model and represent the overall stack heat capacity, lye circulation heat capacity, number of electrolytic cells, thermoneutral voltage, active transfer area, heat transfer coefficient, reversible voltage, and electrodes area, respectively. $r_1$, $r_2$, $t_1$, $t_2$, $t_3$, and $s$ are the parameters related to the electrochemical model. We use the values for the parameters provided by \cite{Rizwan:etal:2021}. The lye inlet flow rate $u(t)=f_{in}$, is the input. The disturbances, $d(t) = \begin{bmatrix}
 T_{\mathrm{amb}};& P_{\mathrm{in}}
\end{bmatrix}$, are the ambient temperature and the inlet power, respectively, with the values $T_{\mathrm{amb}} = 25$ $C^{\circ}$ and $P_{\mathrm{in}} = 2$ MW. We assume these disturbances to be known by the NMPC.

We measure the temperature, i.e.,  $y_k^m = T_k+v_k$ in  \eqref{subeq:meas_eq}, where $v_k$ is the measurement noise with covariance $R_k = 1$. We sample the measurements at every $T_s = 4$ min.  The NMPC tracks a time-varying setpoint for the temperature, $T$. Hence, The output equation in \eqref{subeq:out_eq} is $z(t) = T$. For the OCP, we apply the control horizon $T_N = 100$ min using $N = 25$ steps. The weight matrices for tracking error and input-rate-of-movement penalties are $Q_z = 10$ and $Q_{\Delta u} = 0.1$, respectively. We constraint the inlet flow to $2 \;\text{kg/s} \leq f_{in} \leq 10 \;\text{kg/s}$ and apply the ESDIRK34 with fixed step size $h_k = h = 0.2T_s$ in both the OCP and CD-EKF. We use the implicit-explicit method to perform stochastic simulations of the closed-loop system.

\begin{figure}[tb]
    \centering
    \includegraphics[width=0.5\textwidth, trim={0.55cm 0.6cm 0.8cm 0.8cm},clip]{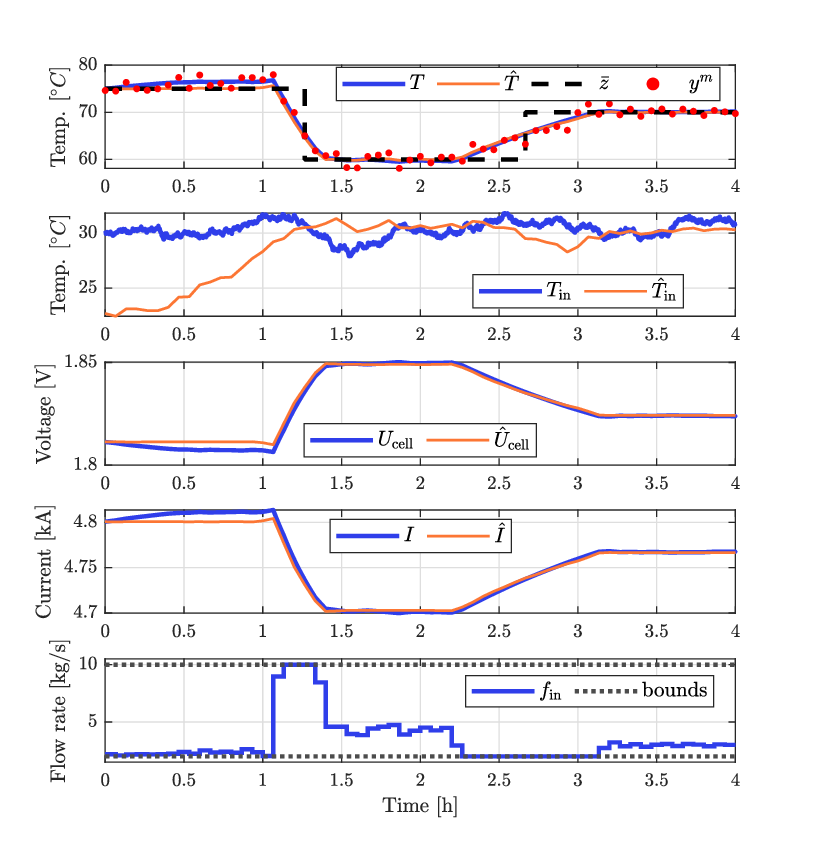}
    \caption{Closed-loop simulation of the alkaline electrolyzer model with NMPC based on ESDIRK34 integration.}
    \label{fig:setpoint_tracking}
\end{figure}

Fig \ref{fig:setpoint_tracking} shows the closed-loop simulation of the alkaline electrolyzer stack, tracking a time-varying setpoint $\Bar{z}(t)$ varying between 75$^{\circ}$C and 60$^{\circ}$C. 
The result shows that the NMPC effectively follows the setpoint, even when the initial estimate of the inlet temperature, $\hat{T}_{\mathrm{in}}$, deviates from the actual state. Improved tracking accuracy is observed as the estimation error in $T_{\mathrm{in}}$ diminishes.

\section{Conclusions}
\label{sec:conclusions}

\noindent
We present an NMPC algorithm for setpoint tracking in systems described by continuous-discrete semi-explicit SDAEs of index 1. The NMPC consists of a CD-EKF and an OCP discretized by direct multiple shooting. We apply ESDIRK integration schemes to solve the IVPs in the CD-EKF and OCP. The ESDIRK implementation uses an iterative IND approach to compute sensitivities. These sensitivities are used to provide accurate gradient information to the OCPs and to compute the covariance in the CD-EKF. We apply the NMPC to an alkaline electrolyzer stack model to track a time-varying setpoint for temperature. 

\bibliographystyle{IEEEtran}
\bibliography{ref/references_cdc2024}

\begin{thebibliography}{10}
\providecommand{\url}[1]{#1}
\csname url@samestyle\endcsname
\providecommand{\newblock}{\relax}
\providecommand{\bibinfo}[2]{#2}
\providecommand{\BIBentrySTDinterwordspacing}{\spaceskip=0pt\relax}
\providecommand{\BIBentryALTinterwordstretchfactor}{4}
\providecommand{\BIBentryALTinterwordspacing}{\spaceskip=\fontdimen2\font plus
\BIBentryALTinterwordstretchfactor\fontdimen3\font minus \fontdimen4\font\relax}
\providecommand{\BIBforeignlanguage}[2]{{%
\expandafter\ifx\csname l@#1\endcsname\relax
\typeout{** WARNING: IEEEtran.bst: No hyphenation pattern has been}%
\typeout{** loaded for the language `#1'. Using the pattern for}%
\typeout{** the default language instead.}%
\else
\language=\csname l@#1\endcsname
\fi
#2}}
\providecommand{\BIBdecl}{\relax}
\BIBdecl

\bibitem{Ritschel:Joergensen:2019}
T.~K.~S. Ritschel and J.~B. Jørgensen, ``Nonlinear {M}odel {P}redictive {C}ontrol for {D}isturbance {R}ejection in {I}soenergetic-isochoric {F}lash {P}rocesses,'' in \emph{IFACPapersOnLine}, vol.~52, no.~1, 2019, pp. 796--801.

\bibitem{Biegler:etal:2012}
L.~T. Biegler, S.~L. Campbell, and V.~Mehrmann, \emph{Control and {O}ptimization with {D}ifferential-{A}lgebraic {C}onstraints}.\hskip 1em plus 0.5em minus 0.4em\relax SIAM: Society for Industrial and Applied Mathematics, 2012.

\bibitem{Ascher:Petzold:1998}
U.~M. Ascher and L.~R. Petzold, \emph{Computer {M}ethods for {O}rdinary {D}ifferential {E}quations and {D}ifferential-{A}lgebraic {E}quations}.\hskip 1em plus 0.5em minus 0.4em\relax SIAM: Society for Industrial and Applied Mathematics, Jan. 1998.

\bibitem{Kristensen:etal:2004}
M.~R. Kristensen, J.~B. Jørgensen, P.~G. Thomsen, and S.~B. Jørgensen, ``An {ESDIRK} method with sensitivity analysis capabilities,'' \emph{Computers \& Chemical Engineering}, vol.~28, no.~12, pp. 2695--2707, 2004.

\bibitem{Bock:1981}
H.~G. Bock, ``{N}umerical {T}reatment of {I}nverse {P}roblems in {C}hemical {R}eaction {K}inetics,'' in \emph{Modelling of Chemical Reaction Systems}, K.~H. Ebert, P.~Deuflhard, and W.~J{\"a}ger, Eds.\hskip 1em plus 0.5em minus 0.4em\relax Berlin, Heidelberg: Springer Berlin Heidelberg, 1981, pp. 102--125.

\bibitem{Albersmeyer:2010}
J.~Albersmeyer, ``Adjoint based algorithms and numerical methods for sensitivity generation and optimal control of large scale dynamic systems,'' Ph.D. dissertation, Ruprecht-Karls-Universität Heidelberg, Dec. 2010.

\bibitem{Christensen:Joergensen:2024}
A.~H.~D. Christensen and J.~B. Jørgensen, ``A {C}omparative {S}tudy of {S}ensitivity {C}omputations in {ESDIRK}-{B}ased {O}ptimal {C}ontrol {P}roblems,'' in \emph{Accepted for the 22nd European Control Conference (ECC24)}, 2024.

\bibitem{Joergensen:ecce:2007}
J.~B. Jørgensen, ``A critical discussion of the continuous-discrete extended {K}alman filter,'' in \emph{European Congress of Chemical Engineering-6}, Copenhagen, Denmark, 2007.

\bibitem{Joergensen:inbook:2007}
J.~B. Jørgensen, M.~R. Kristensen, P.~G. Thomsen, and H.~Madsen, \emph{Assessment and {F}uture {D}irections of {N}onlinear {M}odel {P}redictive {C}ontrol}, ser. Lecture Notes in Control and Information Sciences.\hskip 1em plus 0.5em minus 0.4em\relax Springer Berlin, Heidelberg, 2007, vol. 358, ch. New extended {K}alman filter algorithms for stochastic differential algebraic equations, pp. 359--366.

\bibitem{Kristensen:etal:2005}
M.~R. Kristensen, J.~B. Jørgensen, P.~G. Thomsen, M.~L. Michelsen, and S.~B. Jørgensen, ``Sensitivity analysis in index-1 differential algebraic equations by {ESDIRK} methods,'' in \emph{16th IFAC World Congress}, vol.~38, no.~1, Prague, Czech Republic, 2005, pp. 212--217, 16th IFAC World Congress.

\bibitem{Joergensen:arxiv:2018}
J.~B. Jørgensen, M.~R. Kristensen, and P.~G. Thomsen, ``A {F}amily of {ESDIRK} {I}ntegration {M}ethods,'' \emph{arXiv preprint arXiv:1803.01613}, 2018.

\bibitem{Higueras:Roldan:2005}
I.~Higueras and T.~Roldán, ``Starting {A}lgorithms for a {C}lass of {RK} {M}ethods for {I}ndex-2 {DAE}s,'' \emph{Computers \& Mathematics with Applications}, vol.~49, pp. 1081--1099, 2005.

\bibitem{Bock:etal:2007}
H.~G. Bock, M.~Diehl, P.~Kühl, E.~Kostina, J.~P. Schlöder, and L.~Wirsching, \emph{Assessment and {F}uture {D}irections of {N}onlinear {M}odel {P}redictive {C}ontrol}, ser. Lecture Notes in Control and Information Sciences.\hskip 1em plus 0.5em minus 0.4em\relax Springer Berlin, Heidelberg, 2007, vol. 358, ch. Numerical {M}ethods for {E}fficient and {F}ast {N}onlinear {M}odel {P}redictive {C}ontrol, pp. 163--179.

\bibitem{Cantisani2023}
N.~Cantisani, J.~Dovits, and J.~B. Jørgensen, ``Dynamic modeling of an alkaline electrolyzer plant for process simulation and optimization,'' \emph{arXiv preprint arXiv:2311.09882}, 2023.

\bibitem{Rizwan:etal:2021}
M.~Rizwan, V.~Alstad, and J.~Jäschke, ``Design considerations for industrial water electrolyzer plants,'' \emph{International Journal of Hydrogen Energy}, vol.~46, pp. 37\,120--37\,136, 2021.

\end{thebibliography}

\end{document}